\input amsppt.sty
\loadbold
\TagsOnRight
\hsize 30pc
\vsize 47pc
\def\nmb#1#2{#2}         
\def\cit#1#2{\ifx#1!\cite{#2}\else#2\fi} 
\def\totoc{}             
\def\ign#1{}             

\redefine\g{\frak g}

\redefine\a{\frak a}
\redefine\W{{\Cal W}}
\redefine\b{\frak b}

\redefine\H{{\bar H}}
\redefine\K{{\bar K}}
\redefine\su{\frak su}
\redefine\so{\frak so}
\redefine\sp{\frak sp}
\redefine\u{\frak u}
\redefine\k{\frak k}
\redefine\h{\frak h}
\redefine\bh{{\overline{\frak h}}}

\redefine\m{\frak m}
\redefine\n{\frak n}
\redefine\p{\frak p}

\redefine\z{\frak z}

\redefine\1{{1\over 2}}

\define\row#1#2#3{#1_{#2},\ldots,#1_{#3}}

\topmatter
\title Totally geodesic orbits of isometries\endtitle
\author F. Podest\`a and L. Verdiani  \endauthor 
\email podesta@udini.math.unifi.it, verdiani@udini.math.unifi.it
\endemail
\address University of Florence, Ist. Matematico and Dip. di Matematica 
"U.Dini", Italy\endaddress 
\thanks The authors heartily thank Prof. D. Alekseevsky for useful 
conversations\endthanks 
 \abstract We study cohomogeneity one Riemannian manifolds and we establish
some simple criterium to test when a singular orbit is totally geodesic. 
As an application, we classify compact, positively curved Riemannian
manifolds which are acted on isometrically by a non semisimple Lie group
with an hypersurface orbit.\endabstract

\leftheadtext{\smc }
\rightheadtext{\smc }

\endtopmatter
\document
\subhead\totoc\nmb0{1}.  Introduction \endsubhead
\bigskip
We deal with complete Riemannian manifolds ($M,g$) which are acted on isometrically
by a Lie group $G$, which is closed in the full isometry group of ($M,g$) and 
which has a hypersurface orbit; we will say that the action of $G$ is of 
{\it cohomogeneity one\/}. Such $G$-manifolds have been investigated by several
authors (see e.g.[AA],[AA1] and [AP] for a large bibliography); their
large degree of symmetry allows to construct several interesting examples
of geometric structures which are rare in the homogeneous case, like
Kaehler-Einstein metrics, Riemannian metrics with exceptional holonomy.\par
For general results about cohomogeneity one Riemannian manifolds, we refer the
reader to the basic papers [AA] and [AA1]. We summarize here some basic facts.\par
Given a $G$-manifold ($M,g$) of cohomogeneity one, then the orbit space $M/G$ is
a one dimensional manifold, possibly with boundary; the boundary points
corresponds to singular orbits, which can be one or two. In case $M$ is 
compact, either all orbits are regular or we have exactly two singular orbits 
$B_1,B_2$; if we fix a normal geodesic $\gamma$, i.e. a geodesic $\gamma:\Bbb
R\to M$ which intersects every orbit orthogonally, we may denote by $K$ the 
common stabilizer of all regular points of $\gamma$ and by $H,H'$ the two 
singular stabilizer, so that $K\subset H\cap H'$ and $B_1=G/H,B_2=G/H'$. The 
slice theorem ([Br]) asserts in this case that $H/K$ and $H'/K$ are 
diffeomorphic to spheres. We
will call $\theta$ the triple of subgroups $(H,K,H')$ which we have constructed 
using $\gamma$; it is easy to see that, if we change the choice of the geodesic,
the triple changes under conjugation by some element of $G$. A triple $\theta=
(H,K,H')$ of subgroups of $G$ will be called admissible if $K\subset H\cap H'$
and $H/K,H'/K$ are diffeomorphic to spheres.  One of the main results in [AA]
states that there is a one to one correspondence between $G$-manifolds of
cohomogeneity one with two singular orbits (up to $G$-diffeomorphism)  and
admissible triples (up to conjugation).\par
Our first aim is to give some simple criterium to test when a singular orbit 
is totally geodesic; this is obtained in \S 2, where in Theorem 2.1 we find that 
a bound on the dimension of the singular orbit forces it to be totally geodesic.
Another useful criterium is given essentialy in Lemma 2.3, where we state that
the singular orbit turns out to be totally geodesic as soon as the $G$-action on
it has some non discrete kernel. As a first application of these results, we
prove in Theorem 2.5 that a singular orbit which is, as a homogeneous space, a 
symmetric space of rank one is automatically totally geodesic, with the only
exceptions of the two-dimensional sphere and the projective planes over $\Bbb
{C,H,C}a$.\par
In \S 3, we apply the basic idea stated in Lemma 2.3 to the case of a compact 
positively curved Riemannian manifold which is of cohomogeneity one with respect
to a compact, non semisimple Lie group of isometries. \par
Compact Riemannian manifolds of positive curvature with a non discrete isometry
group have been investigated by Hsiang and Kleiner ([HK]) in the four dimensional
case; in [GS], Grove and Searle investigated the case of the action of some
torus on positively curved manifolds and they established a classification
theorem, up to diffeomorphism, when there exists a $T^1$-action or an
$SU(2)$-action with fixed point set of codimension 2 or not less than $4$
respectively. The problem of classifying compact, positively curved Riemannian
manifolds has been solved in dimension 5 and 6 by Searle ([S]), giving a
complete classification up to diffeomorphism. \par
Our main result is stated in Theorem 3.1, where we classify, for $G$ compact non
semisimple, all possible $G$-manifolds of cohomogeneity one which can carry an 
invariant Riemannian metric with positive curvature. We find that $M$ 
must be diffeomorphic to a compact rank one symmetric space, with the 
exception of the Cayley projective plane.
\bigskip
\subhead\totoc\nmb0{2}. Totally geodesic orbits \endsubhead
\bigskip
Throughout the following, ($M,g$) will be a (complete) Riemannian manifold which is
acted on isometrically by a connected Lie group $G$; the group $G$ will be always
supposed to be closed in the full group of isometries of ($M,g$).  We will
denote by $M_{\text{reg}}$ the open dense subset of $M$ given by regular points
for the $G$-action. Our
first result on the geometry of singular orbits is as follows.
\bigskip
\proclaim {Theorem 2.1} Let ($M,g$) be a Riemannian $G$-manifold, where $G$
is a connected Lie group acting (almost) effectively and isometrically on ($M,g$)
with cohomogeneity $1$. Let $G/H$ be an isolated singular orbit and denote by
$k$ the cohomogeneity of the $H$-action on $G/H$. If
$$\dim G/HÊ< {1\over 2} [\dim M + k - 1],$$
then $G/H$ is totally geodesic.\endproclaim
\bigskip
\demo{Proof} We denote by $K$ a regular isotropy subgroup and we fix a singular
point $q\in M$ such that the isotropy $G_q=H$ contains $K$. We now suppose that
$G/H$ is not totally geodesic, so that there exists a tangent vector $v\in
T_q(G/H)$ such that the geodesic $\gamma:[0,\epsilon)\to M$ (for 
some $\epsilon\in \Bbb R^+$) starting from $q$ with
initial vector $v$ is not entirely contained in $G/H$. Since $G/H$ is supposed
to be an isolated singular orbit, it follows that $\gamma((0,\epsilon))$ intersects
non trivially $M_{\text{reg}}$. We now denote by $\imath$ the isotropy
representation of $H$ and consider the stabilizer $H_v=\{h\in H;
\imath(h)v=v\}.$ It is clear that every element of $H_v$ fixes $\gamma$
pointwise, hence $H_v$ is contained in some conjugate of $K$. Hence $\dim
H_v\leq \dim K$. We therefore get  $$\dim G/H - \dim H + \dim H_v \leq \dim G/H
- \dim H + \dim K. $$ If we denote by $k$ the cohomogeneity of the $H$-action on
$G/H$, this is equal to the cohomogeneity of the $H$-action on the tangent space
$T_q(G/H)$ via isotropy representation and therefore  
$$k\leq \dim G/H - \dim H/H_v \leq \dim G/H - \dim H/K.$$ 
We now recall that the cohomogeneity of the $G$-action on $M$ is equal to the
cohomogeneity of the slice representation of $H$; since $H/K$ is a regular
orbit for the slice representation of $H$, we have that 
$$\dim H/K = \dim M - \dim G/H - p, $$
hence
$$ k \leq  2 \dim G/H - \dim M + p$$
which contradicts our hypothesis.\qed\enddemo
\bigskip
When the action of the group $G$ has cohomogeneity one, the existence of totally geodesic 
orbit can be deduced by purely algebraic and topological assumptions, as the following proposition 
shows.
\bigskip
\proclaim {Proposition 2.2} Let $M$ be a compact manifold of positive Euler
characteristic. If a  compact, non semisimple Lie group $G$ acts almost
effectively on $M$ with cohomogeneity one, then there exists at least one singular
orbit which is totally geodesic with respect to any $G$-invariant Riemannian
metric on $M$.\endproclaim \medskip
\demo{Proof} If the Euler characteristic $\chi(M)$ is positive, then the orbit space $M/G$ 
is homeomorphic to a closed interval $[0,1]$ and there are exactly two singular orbits $B_1,B_2$;
moreover we have that 
$$\chi(M) = \chi(B_1) + \chi(B_2),$$
since a regular orbit is odd-dimensional (see e.g. [AP]). It then follows that at least one singular
orbit, say $B_1=G/H$, has positive characteristic; hence the stability subgroup $H$ is of maximal
rank in $G$. If $G$ is not semisimple, it has some connected center $T$ of positive dimension, which
is contained in $H$; therefore, the action of $G$ on $G/H$ has some non trivial kernel. The
following useful Lemma will conclude the proof.\enddemo
\bigskip
\proclaim {Lemma 2.3} If $G/H$ is not totally geodesic, then the action of $G$ on 
$G/H$ is almost effective (i.e. any ideal of $G$ which is contained in $H$ is finite)
\endproclaim
\medskip
\demo{Proof} We suppose that the representation $\imath$ of $H$ has some 
kernel $N$. If $G/H$ is not totally geodesic, then there exists a tangent 
vector $v\in T_p(G/H)$ such that the subgroup $H_v=\{h\in H; \imath(h)v=v\}$ is 
contained in some conjugate $gKg^{-1}$ for some $g\in G$. But $N\subset H_v$, 
hence $N=g^{-1}Ng\subset K$. If we denote by $\nu$ the slice 
representation of $H$, we have that $\nu(N)\subset \nu(K)$ is 
a normal subgroup of 
$\nu(H)$; since $\nu(H)$ acts effectively on the unit sphere of the normal space 
$T_p(G/H)^\perp$, it follows that $N\subset \ker \nu$. But the action of $G$ on $M$ was supposed to
be almost effective, so  $N=N\cap \ker\nu$ is a finite subgroup of $G$.\qed\enddemo
\bigskip
Again, when the cohomogeneity of the $G$-action is $1$ with one singular orbit $G/H$ and the
cohomogeneity $k$ of the action of $H$ on $G/H$ is $1$, that is when the singular orbit is a
two point-homogeneous space, the estimate given in Corollary 2.2 can be refined. Actually, we may
prove the following \bigskip 
\proclaim{Theorem 2.4} Let ($M,g$) be a complete Riemannian manifold which
is acted on almost effectively and isometrically by a connected Lie group $G$ with
cohomogeneity one. If a singular orbit $G/H$ is a two-point-homogeneous space ,
then $G/H$ is totally geodesic in $M$ unless $G/H$ is the two dimensional sphere or a 
projective space of dimension $2$ over $\Bbb {R,C,H,C}a$.\par
In Table 1 we indicate the exceptional homogeneous
space $G/H$ for compact $G$, a compact Riemannian $G$-manifold $M$ of cohomogeneity one having $G/H$
as non totally geodesic singular orbit and the slice representation $\nu$ of $H$.\endproclaim 
\medskip
\remark {Remark} The manifold $M$ is not necessarily unique; but if ($N,h$) is any Riemannian
manifold  which is acted on by $G$ with non totally geodesic singular orbit $G/H$, then a
$G$-invariant tubular neighborhood of $G/H$ in $N$ is $G$-diffeomorphic to a $G$-invariant
tubular neighborhood of $G/H$ in $M$. The embeddings of the projective spaces $G/H$ into $M$ are the
Veronese embeddings \endremark    
$$\vbox{\offinterlineskip 
\halign {\strut\vrule \hfil \ $#$\ \hfil &\vrule\hfil \  $#$\ \hfil 
&\vrule \hfil\ $#$\ \hfil 
&\vrule\hfil\ $#$\ \hfil&\vrule\hfil\ $#$\ \hfil\vrule\cr 
\noalign{\hrule } 
\ {G/H}_{\phantom {\sum_1^N}} &  
 \ G_{\phantom {\sum_1^N}} &  
 \ H_{\phantom {\sum_1^N}} 
 & \ {M}_{\phantom {\sum_1^N}}^{\phantom 
{\sum_1^N}}& \nu \cr 
\noalign{\hrule depth 1 pt} 
\ {S^2\cong {\Bbb {CP}}^1} &
 {SO(3),SU(2)}
 & \ {SO(2),U(1)} &\  {{\Bbb {CP}}^2} & {\nu(e^{i\theta})=e^{2i\theta}}\cr 
\noalign{\hrule}
\ {{\Bbb {RP}}^2} &
\ {SO(3)}
 & \ {S(O(1)\times O(2))} &\  {S^4}
 & {\nu:O(2)\to O(2)/{\Bbb Z}_2}\cr 
\noalign{\hrule}
\ {{\Bbb {CP}}^2} &
 \ {SU(3)}
 & \ {U(2)} &\  {S^7} & {\nu:U(2)\to SU(2)/{\Bbb Z}_2}\cr 
\noalign{\hrule}
\ {{\Bbb {HP}}^2} &
\ {Sp(3)}
 & \ {Sp(1)\times Sp(2)} &\  {S^{13}}
 & {\nu:Sp(1)\cdot Sp(2)\to SO(5)}\cr 
\noalign{\hrule}
\ {{\Bbb C}a{\Bbb P}^2} &
\ {F_4}
 & \ {Spin(9)} &\  {S^{25}} & {\nu:Spin(9)\to SO(9)}\cr 
\noalign{\hrule} }}$$
\botcaption{Table 1} \endcaption
\demo{Proof} We start recalling that if $p$ is a singular point with singular 
isotropy subgroup $H$,
then the isotropy representation of $H$ at $p$ splits into the sum of two
representations $\nu$ and $\tau$, where $\nu$ is the slice representation 
and $\tau$ is the isotropy representation of $H$ on the homogeneous space
$G/H$. The singular orbit $G/H$ is supposed to be a two-point-homogeneous space, on
which $G$ can be supposed to act almost effectively by Lemma 2.4, if $G/H$ is not totally
geodesic.\par
So, by the classification of two-point-homogeneous spaces, $G$ is a simple Lie group
and  $G/H$ is a symmetric space of rank one. We will analyze the case when $G$ is
compact, while the not compact case follows by duality. Now, for each compact 
symmetric space of rank one (CROSS) $G/H$, we look at all possible
representations $\nu:H\to O(V)$ such that $\nu(H)$ acts transitively and
effectively on the unit sphere of $V$. Of course, we always have the choice
$\nu=\tau$.
\medskip
\proclaim {Lemma 2.6} If $\nu=\tau$, then $G/H$ is totally geodesic.
\endproclaim 
\medskip
\demo{Proof} Indeed, we denote by $h$ the second fundamental form of $G/H$ at
$p\in M$; we may consider $h$ as an element of $S^2(V^{*})\otimes V$, which is
acted on by $H$ in a natural way by $S^{2}\tau^{*}\otimes\tau$. Now, it easy to 
see that, for all $H$
appearing as isotropy subgroups of any CROSS, the space $S^2(V^{*})\otimes V$ has
no non trivial $H$-fixed vectors, so that $h=0$.\qed \enddemo
\bigskip
We may now consider each CROSS separately; for each of them, we look for 
representations $\nu:H\to O(V)$, which are transitive on the unit sphere of
$V$ and for which the $G$-manifold $G\times_{(H,\nu)}V$ admits a Riemannian 
metric such that $G/H$ is not totally geodesic. Throughout the following, we
will keep the same notations, meaning $K$ for a regular isotropy subgroup and 
$H_v$ for a regular isotropy subgroup of $H$ for the action of $H$ on $G/H$;
this means that, by the proof of the main theorem, $H_v$ must be contained in
some conjugate of $K$. Moreover, we recall that we allow the action of $G$ on
$G\times_{(H,\nu)}V$ to be almost effective.  \par  (1)\ $G/H=SO(n+1)/SO(n) =
S^n$. In this case, we may suppose that $G=Spin(n+1)$ and $H=Spin(n)$ for $n\geq
3$ and $H=SO(2)$ for $n=2$.\par Then, any representation $\nu:H\to O(V)$  which
is transitive on the unit sphere of $V$, factors through $\pi:Spin(n)\to SO(n)$
giving rise to the standard representation $\nu=\tau$ of $SO(n)$, unless
$n=2,4,5,6,7,9$. We deal with each case separately.\par 
\noindent (a)\ In case $n=2$, we
have the covering homomorphisms $\nu_k:SO(2)\to SO(2)$ given by the $k$-power,
that is $\nu(e^{i\theta})=e^{ik\theta}$. Again, if we denote by $h$ the
fundamental form at $p$ and if we denote by $\sigma$ the element of $H=SO(2)$
such that $\nu_k(\sigma)=-Id$, then, for all tangent vectors $X,Y\in T_p(G/H)$,
we have that $h(X,Y) = - h(\tau(\sigma)X,\tau(\sigma)Y)$. Now, it is
elementary to see that if $k\not=2$, then $h=0$. We now construct an example of a
$G=SO(3)$-action on some Riemannian manifold $M$ such that one  singular orbit
$G/H$ is $S^2$ and is not totally geodesic. In order to do this, we consider the 
action of $SO(3)$ on
${\Bbb CP}^2$ given by the standard inclusion $SO(3)\subset U(2)$; 
then the singular orbit $G/H$ is
the orbit of $G$ through the point $[(1,0,i)]$, which is easily seen to be 
not totally geodesic.\par 
\noindent (b)\ In case $n=4$, we have that $Spin(4)\cong Sp(1)\times Sp(1)$ and 
we may take $\nu:Spin(4)\to Sp(1)$, with $V={\Bbb H}$. In this case, $K$ can
be taken to be one normal factor $Sp(1)$ of $Spin(4)$, while $H_v$ can be taken to be a copy
of $Sp(1)$ embedded diagonally into $Spin(4)$. Now, $H_v$ is not contained in
any conjugate $gKg^{-1}$ for $g\in Spin(5)$, as one can easily see by considering
the projected subgroups $\pi(K)$ and $\pi(H_v)$ of $SO(5)$, where
$\pi:Spin(5)\to SO(5)$ is the standard covering map:\par 
\noindent (c)\ When $n=5,6$,
we have special isomorphisms, namely $Spin(5)\cong Sp(2)$
and $Spin(6)\cong SU(4)$, so that we can choose $V={\Bbb H}^2$ and 
$V={\Bbb C}^4$ respectively. Here Theorem 2.1 applies and the singular orbit is
totally geodesic. When $n=7,9$, we have special representations of $Spin(7)$ 
and $Spin(9)$ which are transitive on the unit sphere of $\Bbb R^8$ and 
$\Bbb R^{16}$ resp. Again Theorem 2.1 applies.\par (2)\
$G/H=SU(n+1)/S(U(1)\times U(n))= {\Bbb CP}^n$. Here we have $G=SU(n+1)$ and
$H=U(n)$ and we have several possibilities for the representation $\nu$. \par 
\noindent (a)\ $\nu:U(n)\to T^1/{\Bbb Z}_n\cong T^1$. In this case $K=\ker
\nu=SU(n)$. If the orbit $G/H$ should not be totally geodesic, then a regular
isotropy subgroup $H_v$ for the action of $H$ on $G/H$ should be contained in
some conjugate of $K$; but this is easily seen to be impossible.\par 
\noindent (b)\ We have
$\nu_k:U(n)\to U(n)$ given by $\nu_k(A)=(\det A)^kA$ for $k\in  \Bbb Z$; they
represent, up to equivalence, all possible coverings of $U(n)$. But again an
easy computation shows that the subgroup $H_v$ is contained in some conjugate of
$K$ if and only if $k=1$; so by Lemma 2.4, the orbit $G/H$ is totally geodesic.
\par 
\noindent (c)\ When $n=2$, we have $\nu:U(2)\to SU(2)/{\Bbb Z}_2\cong SO(3)$. In
this case, in order to construct the manifold $M$, it is enough to consider the
sphere $S^7\subset \su(3)$, which is acted on by $SU(3)$ by adjoint
representation; this action has cohomogeneity one and one singular orbit is
${\Bbb CP}^2$, which cannot be totally geodesic in $S^7$.\par 
(3)\ $G/H=Sp(n+1)/Sp(1)\cdot Sp(n)\cong {\Bbb HP}^n$. Here we may take the simply
connected group $G=Sp(n+1)$ and $H=Sp(1)\times Sp(n)$. We have to consider the
following possibilities:\par  
\noindent (a)\ $\nu:Sp(1)\times Sp(n)\to
Sp(1)/{\Bbb Z}_2=SO(3)$. Here $K=T^1\times Sp(n)$ and $H_v=Sp(1)\times 
Sp(n-1)$, which is not contained in any conjugate of $K$. So in this case the 
orbit $G/H$ is totally geodesic. \par 
\noindent (b)\ $\nu:Sp(1)\times Sp(n)\to Sp(n)$; in this case $H_v=Sp(1)\times
Sp(n-1)$ and $K=Sp(1)^{\Delta}\times Sp(n-1)$, where $Sp(1)^{\Delta}$ is a copy
of $Sp(1)$ embedded diagonally into $Sp(1)\times Sp(1)\subset Sp(1)\times
Sp(n)$. Here again $H_v$ is not contained in any conjugate of $K$, since, as
subgroups of $Sp(n+1)$, we have that Fix(${\Bbb H}^{n+1},H_v$) is a
quaternionic line, while Fix(${\Bbb H}^{n+1},K$) is reduced to $\{0\}$.\par
\noindent (c)\ $\nu:Sp(1)\times Sp(2)\to Sp(2)/{\Bbb Z}_2\cong SO(5)$.
In this case the manifold $M$ can be chosen to be the sphere $S^{13}$; indeed the
group $Sp(3)$ admits a cohomogeneity two linear action given by
$\Lambda^2\nu_3-1$, with representation space ${\Bbb R}^{14}$ (see[HL]); one singular
orbit for the action on the unit sphere is a copy of the quaternionic projective
plane ${\Bbb HP}^2$ (Veronese embedding).\par 
(4)\ $G/H=F_4/Spin(9)={\Bbb C}a{\Bbb P}^2$. In this
case the only representation $\nu$ different from $\tau$ is given by
$\nu:Spin(9)\to SO(9)$. If we now consider the irreducible linear
representation $\phi$ of $F_4$ of dimension $26$, it has cohomogeneity two (see
[HL]) and one singular orbit for the $F_4$-action on $S^{25}$ is a copy of a
Cayley projective plane (again Veronese embedding).\qed
\enddemo
\bigskip
\subhead\totoc\nmb0{3}.  Manifolds of positive curvature \endsubhead
\bigskip
In this section, we will use some ideas which were developped in the 
previous section and we will consider even dimensional, compact Riemannian
manifolds ($M,g$) of positive curvature, which are acted 
isometrically and (almost) effectively on by a compact, non
semisimple Lie group $G$ of isometries with one hypersurface orbit. \par
Our aim is to prove the following
\proclaim{Theorem 3.1} Let $M^{2n}$ be an even dimensional, compact 
Riemannian manifold of positive sectional curvature. If a compact, non 
semisimple Lie 
group $G$ acts isometrically on $M^{2n}$ by
cohomogeneity one, then \roster
\item the dimension of the center of $G$ is one and the semisimple 
part of $G$ acts by cohomogeneity one;
\item the manifold $M$ is diffeomorphic to a rank one symmetric space 
with the exception of the Cayley projective plane.\endroster
\endproclaim\medskip
We start proving some basic lemmata.
\bigskip
\proclaim {Lemma 3.2} Let $M^{2n}$ be an even dimensional, compact 
Riemannian manifold of positive sectional curvature. If a compact Lie 
group $G$ acts isometrically on $M^{2n}$ by
cohomogeneity one, then $M^{2n}$ has positive Euler characteristic and 
there are exactly two singular orbits, one of which has positive Euler 
characteristic.\endproclaim
\bigskip
\demo{Proof} We note that $M^{2n}$ has finite fundamental group, so 
that there is no fibration $M^{2n}\to S^1$. It then follows that there 
are exactly two singular orbits, $B_1,B_{2}$. We now claim that at 
least one of these has positive Euler characteristic, which will 
conclude the proof. In fact, let $T$ be a maximal torus of $G$ and 
let $X$ be a Killing vector field on $M^{2n}$ so that 
$\overline{\{\exp(tX)\}}=T$. By Berger's Theorem (see e.g.[Ko]), the field 
$X$ has at least one zero $p\in M^{2n}$. We claim that we can suppose 
that $p$ is a singular point. Indeed, if $p$ is regular, then $X$ 
vanishes along the whole normal geodesic issuing from $p$, hence it 
vanishes at some singular point too. Therefore, we have that $T\subset 
G_{p}$ and $G_{p}$ has maximal rank in $G$. This implies that 
$\chi(G/G_p)>0$. \qed \enddemo
\bigskip
\proclaim {Lemma 3.3} Let $M^{2n}$ be an even dimensional, compact 
Riemannian manifold of positive sectional curvature. If a compact, 
non semisimple Lie  group $G$ acts isometrically and almost effectively on $M^{2n}$
by cohomogeneity one, then \roster
\item the dimension of the center $Z(G)$ is one;
\item there is a singular orbit which is totally geodesic and of positive 
Euler characteristic.
\endroster \endproclaim
\demo {Proof} By Lemma 3.2, we know that there exists a singular orbit $B$ 
which is of positive Euler characteristic; hence, if we represent $B=G/H$,
then $Z(G)\subset H$ and $B$ is totally geodesic by Lemma 2.3. Now $Z(G)$ acts 
trivially on $B$ and, if $\nu$ denotes the slice represenation of $H$, then 
$\nu(Z(G))$ is at most one-dimensional; then almost effectiveness of the
$G$-action implies that $\dim Z(G)=1$.\qed\enddemo
\bigskip
Throughout the following, we will always denote by $B$ the totally
geodesic singular orbit with $\chi(B)>0$. Moreover we will use the following 
results by Grove-Searle ([GS])
\bigskip
\proclaim {Theorem [GS]} Let $M^n$ be a simply connected, compact positively
curved manifold. \roster
\item If $T^1$ acts isometrically and effectively on $M^n$ with a fixed point set 
of codimension $2$, then $M^n$ is diffeomorphic to $S^n$ or to 
${\Bbb {CP}}^{n/2}$;
\item If $SU(2)$ acts isometrically and almost effectively on $M^n$ with 
fixed point set of codimension less or equal to $4$, then $M^n$ is diffeomorphic 
to $S^n,{\Bbb {CP}}^{n/2}$ or ${\Bbb {HP}}^{n/4}$.\endroster\endproclaim
\bigskip
Since the manifold $(M,g)$ has finite fundamental group, we will always assume 
that $M$ is simply connected. Moreover, when the singular orbit $B$ is reduced 
to a point, i.e. when $G$ has a fixed point, we know form [AP], that the manifold
is diffeomorphic to a CROSS (compact rank one symmetric space), while if the 
codimension of $B$ is $2$, we may apply Theorem [GS], to get that $M$ is
diffeomorphic to $S^n$ or to ${\Bbb {CP}}^{n/2}$. So, in the proof of Theorem , 
we will always assume that $\text{codim} B>2$.\par 
We may always decompose the group $G$ as $G = T^1\cdot G_o\cdot G_1$, where 
$T^1\cdot G_o$ is the kernel of the $G$-action on $B$; accordingly, we have that
$H=T^1\cdot G_o\cdot H_1$, where $G_1$ acts (almost) effectively on $B$ with 
$B=G_1/H_1$. We denote by $\nu$ the slice representation of $H$ and by gothic 
letters the Lie algebras.\par
\proclaim {Case $\g_0\not=0$}\endproclaim
In this case we claim that there is an $SU(2)$-action whose fixed point set has 
codimension $4$ in $M$, so that we can apply Theorem [GS].\par
Since the $G$-action on $M$ is almost effective, we have that $\nu|_{\g_o+\Bbb
R}$ is an isomorphism and since $\nu(\g_o)$ has a non trivial centralizer, we
have $\g_0\cong \su(m)$ or $\sp(m)$. We fix invariant complements $\p$ and $\m$ 
so that 
$$\g=\h+\m,\ \g_1=\h_1+\m,\ \h=\k+\p.$$
At a regular point $y\in M$, we have the identification (as $\k$-modules)
$$T_yM = \Bbb R + \p + \m.$$
Moreover, $\nu(\g_o)$ still acts transitively on the unit normal sphere, so we
can consider the stabilizer $\k_o\subset \g_o$, with $\k_o\cong
\su(m-1),\sp(m-1)$. Now, if $\k_o=\{0\}$, then $\g_o\cong\su(2)$ and $G_o^o$,
which is locally isomorphic to $SU(2)$, fixes $B$ which has codimension $4$. \par
If $\k_o\not=\{0\}$, then $\k_o$ acts trivially on $\m$ and $\p$ splits as 
$\p=\p_o+\p_1$, where $\p_o$ is a one- or three-dimensional trivial summand and 
$\p_1$ is $\k_o$-irreducible. In any case, there is a subalgebra
$\su(2)\subset\k_o$, whose fixed point set in $\p$ has codimension $4$ in $\p$
and we are done. \bigskip
\proclaim {Case $\g_o=\{0\}$}\endproclaim
In this case $G=T^1\cdot G_1$ and $H=T^1\cdot H_1$, where $G_1/H_1$ is even
dimensional and carries a positively curved invariant metric. Moreover, $B$ is
the fixed point set of $Z(G)$, hence it is simply connected by Synge Theorem
([Ko]). We give here the list of all such pairs ($G_1,H_1$) together with the
dimension of the corresponding space $G_1/H_1$; we call this the Wallach's list
([Wal]): $$\vbox{\offinterlineskip  \halign {\strut\vrule \hfil \ $#$\ \hfil
&\vrule\hfil \  $#$\ \hfil  &\vrule \hfil\ $#$\ \hfil 
&\vrule\hfil\ $#$\ \hfil\vrule\cr 
\noalign{\hrule } 
\ {n.}_{\phantom {\sum_1^N}} &  
 \ {G_1}_{\phantom {\sum_1^N}} &  
 \ {H_1}_{\phantom {\sum_1^N}} 
 & \ {\dim}_{\phantom {\sum_1^N}}^{\phantom 
{\sum_1^N}}\cr 
\noalign{\hrule depth 1 pt} 
\ {1} &
 {SU(n+1)}
 & \ {U(n)} &\ {2n} \cr 
\noalign{\hrule}
\ {2} &
 {SU(3)}
 & \ {T^2} &\ {6} \cr 
\noalign{\hrule}
\ {3} &
 {Spin(2n+1)}
 & \ {Spin(2n)} &\ {2n} \cr 
\noalign{\hrule}
\ {4} &
 {Sp(n)}
 & \ {Sp(1)\times Sp(n-1)} &\ {4(n-1)} \cr 
\noalign{\hrule}
\ {5} &
 {Sp(n)}
 & \ {T^1\times Sp(n-1)} &\ {2(2n-1)} \cr 
\noalign{\hrule}
\ {6} &
 {Sp(3)}
 & \ {Sp(1)^3} &\ {12} \cr 
\noalign{\hrule}
\ {7} &
 {{\frak F}_4}
 & \ {Spin(9)} &\ {16} \cr 
\noalign{\hrule}
\ {8} &
 {{\frak F}_4}
 & \ {Spin(8)} &\ {24} \cr 
\noalign{\hrule}
\ {9} &
 {{\frak G}_2}
 & \ {SU(3)} &\ {6} \cr 
\noalign{\hrule}}}$$
$$\text{Wallach's list}$$ 
Since we are supposing that ${\text{codim}}\ B>2$, we have that $\nu|_{H_1}$ is
not trivial, so $\h_1$ must contain some ideal isomorphic to $\su(m)$ or to
$\sp(m)$. It then follows that we may exclude cases (2),(7),(8) and (3) (for
$n\not=2,3$) from Wallach's list. \par
Our next arguments will rely on the study of the second singular orbit $B'=G/H'$,
where $K\subset H\cap H'$; we will denote by $\nu':H'\to O(V')$ the slice
representation of $H'$ on the normal space $V'$. The following two observations
will be useful 
\medskip
{\bf Fact 1.} The Lie algebra $\h'$ cannot be semisimple. Otherwise it would be
contained in $\g_1$, while $\k$ always has a not trivial projection on
$\z(\g)$.\par
{\bf Fact 2.} If $\h'$ has maximal rank, then $\dim V'>\dim B$. Indeed, if $\h'$
has maximal rank, then $B'$ is totally geodesic and, by Frankel Theorem ([Ko]),
$\dim B+\dim B'<\dim M$.\par
Now, for each case in Wallach's list, we compute $\k$ and for each ideal $\n$ of 
$\k$, we look for subalgebras $\h'$ in $\g$ such that $\h'$ contains $\k$ and has
$\n$ as an ideal with ($\h'/\n,\k/\n$) belonging to Borel's list.\par
Using Fact 1 and 2, it is not difficult to exclude case (9) and (3) for $n=2$.\par
We are let with the cases (1),(4),(5),(6). Case (1) is handled in the next 
\proclaim {Lemma 3.4} If $G_1=SU(n)$, then the manifold $M$ is diffeomorphic to 
a sphere, a complex or quaternionic projective space.\endproclaim
\bigskip
\demo {Proof} We have that $G_1=SU(m)$ and $H_1=U(m-1)$ for some $m$. Since 
we are supposing that the codimension of $B$ is bigger than $2$, we have that
$\nu(H)=U(m-1)$ and $\nu(H_1)$ still acts transitively on unit normal sphere. 
Note that $\dim M=4(m-1)$, which can be assumed to be bigger than $4$, by the
results of Hsiang and Kleiner ([HK]). So, the $G_1$-action on $M$ is of
cohomogeneity one. The slice representation $\nu|_{H_i}$ can be assumed to be of
the form 
$$\nu((e^{i\theta},A)) = (\det A)^{j}A,\quad 
(e^{i\theta},A)\in S(U(1)\times U(q)),$$
where $j\in \Bbb Z$. It then follows that we can fix $K$ to be 
$$K = \{\left(\matrix e^{i\theta} &{}&{}\cr {}&e^{i\phi}&{}\cr
{}&{}&A\endmatrix\right)\in SU(q+1); 
e^{i(j+1)\phi)}(\det A)^j = 1\},$$
so that $\k\cong \Bbb R+\su(m-2)$. Now, for any ideal $\n'\subset\k$, we determine
all subalgebras $\h'$ of $\g_1$ such that $\k\subset \h'$ and 
($\h'/\n',\k/\n'$) belongs to Borel's list.\par
We have the following possibilities (here $\n'$ and
$\h'$ are given up to isomorphism): \bigskip $$\vbox{\offinterlineskip 
\halign {\strut\vrule \hfil \ $#$\ \hfil &\vrule\hfil \  $#$\ \hfil 
&\vrule \hfil\ $#$\ \hfil 
&\vrule\hfil\ $#$\ \hfil\vrule\cr 
\noalign{\hrule } 
\ n._{\phantom {\sum_1^N}} &  
 \quad \n'_{\phantom {\sum_1^N}} &  
 \quad {\bh'}_{\phantom {\sum_1^N}} 
 & \ {\dim V'}_{\phantom {\sum_1^N}}^{\phantom 
{\sum_1^N}}\cr 
\noalign{\hrule depth 1 pt}
{1} & 
\ {\{0\}} &
 {\Bbb R+\su(m-1)}
 & \ {2q} \cr 
\noalign{\hrule}
{2} &
\ {\k}
 & \ {\Bbb R^2+\su(m-2)} & {2}\cr 
\noalign{\hrule}
{3} &
\ {\k}
 & \ {\Bbb R+\su(2)+\su(m-2)} & {4}\cr 
\noalign{\hrule}
{4} &
 \ {\Bbb R}
 & \ {\Bbb R+\su(m-1)} &{2(m-1)} \cr 
\noalign{\hrule}
{5} &
 \ {\Bbb R, m=4}
 & \ {\Bbb R+\so(4)} &{4} \cr 
\noalign{\hrule}
{6} &
\ {\su(m-2)}
 & \ {\u(2)+\su(m-2)} &{4}\cr 
\noalign{\hrule}
{7} &
\ {\su(m-2)}
 & \ {\su(2)+\su(m-2)} &{3}\cr 
\noalign{\hrule} }}$$
\bigskip
First of all we must exclude the case (1)-(5) of the previous table. 
In each case (1)-(5), the subalgebra $\h'$ is of maximal rank, hence 
$\chi(B')>0$. On the other hand, we have that $G=T^1\cdot G_1$ and, if 
$\exp(tX)\in T^1$, then $\exp(tX)$ maps $B'$ onto a $G_1$-orbit; since $B'$ is 
singular, we get that $T^1$ maps $B'$ onto itself, hence $G(B')=B'$. Since 
$\chi(B')>0$, we may represent $B'$ as $B'=G/H'$ with $H'$ of 
maximal rank in $G$ so that $B'$ is totally geodesic.
But, it can be checked that $\dim G/H + \dim G/H'
\geq 4(m-1)=\dim M$ and by Frankel's theorem (see [Fr]) the two singular orbits
should intersect, which is not the case. \par 
So only case (5) is admissible. Moreover,
since $H$ is connected and $H/K\cong S^{2m-3}$, with  $2m-3>1$, we have that
$K$ is connected; since now $H'/K\cong S^2$, we get that $\H'$ is 
connected.\par
In order to avoid confusion, we put $\h' = \a + \b,$
where $\a,\b$ are ideals 
of $\h'$ which are isomorphic to $\su(2)$ and $\su(m-2)$ 
respectively.\par
We will now devide the discussion according to $m\geq 4$ or $m=3$.\par
\medskip
\proclaim {Case $m\geq 4$}\endproclaim
In this case $\b\not=\{0\}$ and $\a$ centralizes $\su(q-1)$ in $\g_1$; 
so $\a$ is contained in the 
semisimple part of the centralizer $C_{\g_1}(\su(m-2))$, which is 
isomorphic to $\su(2)$. It then follows that $\a$ is embedded into 
$\g_1$ as the set of all matrices
$$\a = \{\left(\matrix \su(2)& 0\cr
0& 0\endmatrix\right) \in \su(m)\}$$
and therefore
$$\h' = \{\left(\matrix \su(2)& 0\cr
0& \su(m-2)\endmatrix\right) \in \su(m)\}.$$
Since $\K\subset \H\cap\H'$, we get that 
$$\K = \{\left(\matrix e^{i\theta} &{}&{}\cr {}&e^{-i\theta}&{}\cr
{}&{}&A\endmatrix\right)\in SU(m); 
\det A = 1\},$$
which forces $j=-1$.\par
Therefore we have proved that the $G_1$-action on $M$ is associated to 
the triple $\theta = (H,K,H')$, which is uniquely given by 
$$\theta = (S(U(1)\times U(m-1)), S(U(1)\times U(1))\times SU(m-2),SU(2)\times
SU(m-2)).$$
On the other hand, if we consider the standard embedding $SU(m)\subset
U(m)\subset Sp(m)$, then the group $SU(m)$ acts on the space 
$Sp(m)/Sp(1)\cdot Sp(m-1)\cong \Bbb {HP}^{m-1}$ by cohomogeneity one and associated
triple of subgroups $\theta$. According to the main results in [AA], we conclude
that $M$ is diffeomorphic to $\Bbb {HP}^q$.
\medskip
\proclaim {Case $m=3$}\endproclaim
In this case $\h'\cong \su(2)$ and 
$$\k = \{{\text{diag}}(i\theta,ij\theta,-i(j+1)\theta)\in \su(3);\ 
\theta\in\Bbb R\}.$$
Up to conjugation there are exactly two imbeddings of $\su(2)$ into $\su(3)$, 
corresponding to a reducible of irreducible representation $\pi$ of $\su(2)$ on 
$\Bbb C^3$; but in any case $\pi$ has $0$ as a weight, so that we must have 
$j=0$ or $j=-1$. Since the two subalgebras $\k$ corresponding to $j=0$ and $j=-1$
are conjugated by some element which normalizes $\h$, we may suppose that $j=-1$.
In this case, there is only one subalgebra $\h'$ isomorphic to $\su(2)$ and
containing $\k$, so that, by the same argument as in the previous case, the 
manifold is diffeomorphic to ${\Bbb {HP}}^2$.\qed\enddemo
\bigskip
We now proceed excluding the remaining cases.
\proclaim {Lemma 3.5} Case (6) cannot occur. \endproclaim
\bigskip
\demo {Proof} In this case $G=T^1\cdot Sp(3)$, $H=T^1\cdot Sp(1)^3$ with 
$\nu(H)=U(2)$. Using Fact 1 and 2, it is not difficult to see that the only 
possibility for $\h'$ is $\h'\cong \Bbb R + \sp(2)$. We will now consider the 
action of the semisimple part of $G$, which still acts by cohomogeneity one with
associated triple $\theta=(H,K,H')=(Sp(1)^3,Sp(1)^2,Sp(2))$. We may write 
$$\sp(3) = \k + \m_o + \m_1 + \m_2 + \m_3,\tag 3.1$$
where $\m_o$ is a trivial $ad(\k)$-module and $\m_i$ are irreducible, mutually
inequivalent $ad(\k)$-modules. We now fix a non zero vector $v\in
\m_1$ and  consider a normal geodesic $\gamma:\Bbb R\to M$ w.r.t. a positively
curved  $G$-invariant metric $g$ on $M$; we choose $\gamma$ so that it induces the 
triple $\theta$. \par
First of all, we claim that the Killing vector field $X$ induced by $v$ on $M$ 
never vanishes along $\gamma$. In order to prove this, it is enough to check that
$v$ does not belong to $Ad(\W)\h\cap Ad(\W)\h'$, where $\W$ is the generalized
Weyl group generated by the two geodesic symmetries $\sigma,\sigma'$. But 
we have that 
$$Ad(\sigma)\m_i=\m_i,\ i=1,2\ {\text{and}}\  Ad(\sigma')\m_1=\m_2,
Ad(\sigma')\m_2=\m_1.$$
We now consider the smooth function $f(t)=||X||_{\gamma(t)}$ for $t\in\Bbb R$
and we claim that $f$ is a concave positive function, which is not possible.
This will conclude our proof. \par
It will be enough to check that $f"(t)<0$ for all $t$ such that $\gamma(t)$ is 
a regular point. First, we observe that (3.1) gives a decomposition of the 
tangent space to a regular orbit into $K$-irreducible, inequivalent submodules,
so that the shape operator of the regular orbit hypersurface will 
preserve each submodule and will be a multiple of the identity operator on
$\m_1$. Therefore, if we denote by $D$ the Levi-Civita connection of $g$, we 
have that $D_{\gamma(t)'}X$ is a multiple of $X_{\gamma(t)}$; we then have 
$$\eqalign {2R_{X\gamma'X\gamma'} &= 2||D_{\gamma'}X||^2 - {{d^2}\over{dt^2}}
f^2\cr
{}&= 2{{g(D_{\gamma'}X,X)^2}\over{f^2}} - 2(f')^2 - 2ff" = -2ff" > 0,\cr},$$
since $g(D_{\gamma'}X,X) = ff'$. \qed\enddemo 
\proclaim {Lemma 3.6} Case (3) for $n=3$ cannot occur.\endproclaim
\bigskip
\demo{Proof} In this case we have $G=T^{1}\times Spin(7)$ and 
$H=T^{1}\times Spin(6)$, with $\nu(H)=U(4)$. Using Fact 1 and 2, it 
is easy to see that the only possibility for $\h'$ is $\h'\cong \Bbb 
R+\su(4)$; moreover there is only one subalgebra of $\so(7)$ which is 
isomorphic to $\su(4)$ and which contains $\k\cong \su(3)$. So, we have two 
totally geodesic singular orbits $B,B'\cong S^{6}$ in codimension $8$ 
with $H=H'$. We now consider the normal space $V$ to the singular 
orbit $B$ at a point $p\in B$ and we claim that there exists an element $\h\in H$ such that
$\nu(h)$ is the identity, while $\tau(h)=-Id$, where $\tau$ denotes 
the tangent isotropy representation.\par
Indeed, $B$ is the symmetric space $Spin(7)/Spin(6)$ and the 
symmetry $\sigma$ belongs to the center of $Spin(6)$, hence 
$\nu(\sigma)$ acts a scalar multiple of the identity on $V$; then there 
is an element $t\in T^{1}$ such that $\nu(t,\sigma)=Id$ and 
$\tau(t,\sigma)=-Id$ and we are done. This means that there exists a 
totally geodesic submanifold $F$ having $V$ as tangent space. The 
submanifold $F$ is of cohomogeneity one under the action of $H$ with
two fixed points $p,p'$ as singular orbits. If now $\gamma$ denotes a normal 
geodesic issuing from $p$, then the parallel transport along $\gamma$
of the tangent space $T_{p}F=V$ will be tangent to $F$; since 
$T_{p'}F=V'$, where $V'$ is the normal space to $B'$ at $p'$, the 
parallel transport along $\gamma$ will send the tangent space $T_{p}B$ 
onto the tangent space $T_{p'}B'$. This is enough to apply Frankel's 
argument ([Fr]) and get a contradiction. \qed\enddemo
The last case is considered in the next 
\bigskip
\proclaim {Lemma 3.7} Case (5) cannot occur. \endproclaim 
\bigskip
\demo{Proof} Here we have $G=T^{1}\times Sp(n)$ and $H= T^{2}\times 
Sp(n-1)$. We will indicate by $B$ the singular orbit $G/H$ and $B'$ 
the second singular orbit; moreover we will indicate by $L$ the 
distance between the two singular orbits, so that the normal 
geodesic $\gamma(t)$ will intersect $B$ for $t\in 2{\Bbb Z} L$.
Using Fact 1 and 2, it is easy to see that the only 
possibility for $\h'$ is $\h'\cong \Bbb R + \sp(n-2) + \su(2)$, where
the normal space $V'$ to the second singular orbit is of dimension 
$3$. We decompose the Lie algebra $\g$ as sum of $\k$-modules, we 
have exactly four irreducible $\k$-modules $\n_{i},\ i=1,\dots,4$ 
of real dimension $2$ and 
pairwise not equivalent. If we choose $X_{i}\in \n_{i}$
for $i=1,\dots,4$, then the functions $f_{i}(t) = ||X_{i}||_{\gamma(t)}$ must be 
smooth concave where they do not vanish, according to the same 
argument as in Lemma 3.5. Moreover three of them do not vanish 
for $t\in 2\Bbb Z L$. Since the geodesic symmetry $\sigma$ at $B$
preserves all these modules, the functions $f_{i}(t)$ which do not vanish
for $t=2\Bbb Z L$, must be even at such values of $t$. Using the 
fact that $\dim V' = 3$, it is now 
easy to see that at least one $f_{i}$ never vanishes, giving a 
contradiction. \qed\enddemo

\enddocument\bye